\documentclass[11pt,epsf]{article}
\usepackage{alltt}
\usepackage[all]{xy}
\usepackage{latexsym,amssymb,amsmath,amsfonts, amsthm, xcolor}

\usepackage{enumerate}
\usepackage{graphicx}
\graphicspath{ {./images/} }

\newtheorem{theorem}{Theorem}
\newtheorem{itlemma}{Lemma}[section]
\newtheorem{itproposition}[itlemma]{Proposition}
\newtheorem{itcorollary}[itlemma]{Corollary}
\newtheorem{itremark}[itlemma]{Remark}
\newtheorem{itremarks}[itlemma]{Remarks}
\newtheorem{itdefinition}[itlemma]{Definition}
\newtheorem{itexample}[itlemma]{Example}

\newenvironment{lemma}{\begin{itlemma}\rm}{\end{itlemma}} 
\newenvironment{remark}{\begin{itremark}\rm}{\end{itremark}} 
\newenvironment{remarks}{\begin{itremarks} \rm}{\end{itremarks}}
\newenvironment{corollary}{\begin{itcorollary}\rm}{\end{itcorollary}}
\newenvironment{proposition}{\begin{itproposition}\rm}{\end{itproposition}}
\newenvironment{definition}{\begin{itdefinition}\rm}{\end{itdefinition}}
\newenvironment{example}{\begin{itexample}\rm}{\end{itexample}}
\newenvironment{fact}{\noindent {{\bf Fact}}:\ \ }{\hfill \medskip}
\newenvironment{claim}{\noindent {\em Claim}. \ \ }{\hfill \medskip}
\newcommand{\be}[1]{\begin{equation}\label{#1}}
\newcommand{\ee}{\end{equation}}
\newcommand{\bl}[1]{\begin{lemma}\label{#1}}
\newcommand{\br}[1]{\begin{remark}\label{#1}}
\newcommand{\brs}[1]{\begin{remarks}\label{#1}}
\newcommand{\bt}[1]{\begin{theorem}\label{#1}}
\newcommand{\bd}[1]{\begin{definition}\label{#1}}
\newcommand{\bp}[1]{\begin{proposition}\label{#1}}
\newcommand{\bc}[1]{\begin{corollary}\label{#1}}
\newcommand{\bfact}[1]{\begin{fact}\label{#1}}
\newcommand{\bex}[1]{\begin{example}\label{#1}}
\newcommand{\ec}{\end{corollary}}
\newcommand{\efact}{\end{fact}}
\newcommand{\eex}{\end{example}}
\newcommand{\el}{\end{lemma}}
\newcommand{\er}{\end{remark}}
\newcommand{\ers}{\end{remarks}}
\newcommand{\et}{\end{theorem}}
\newcommand{\ed}{\end{definition}}
\newcommand{\ep}{\end{proposition}}
\newcommand{\epr}{\end{proof}}
\newcommand{\bpr}{\begin{proof}}
\newcommand{\bcl}{\begin{claim}}
\newcommand{\ecl}{\end{claim}}

\newcommand{\bi}{\begin{itemize}}
\newcommand{\ei}{\end{itemize}}
\newcommand{\ben}{\begin{enumerate}}
\newcommand{\een}{\end{enumerate}}

\title{\bf \Large{On $K-P$ sub-Riemannian Problems and their Cut Locus$\dagger$\thanks{$\dagger$Another shorter version of this paper to appear in Proceedings of the European Control Conference 2019}}
}

\author{Domenico D'Alessandro$^{1}$  \,  and \, Benjamin Sheller$^{2}$
\thanks{$^{1}$ Domenico D'Alessandro is with the Department of Mathematics, Iowa State University, Ames, Iowa, U.S.A., {\tt\small daless@iastate.edu}}
\thanks{$^{2}$ Benjamin Sheller is with Department of Mathematics, Iowa State University, Ames, Iowa, U.S.A., {\tt\small bsheller@iastate.edu}}
}
\begin{document}

\maketitle
\thispagestyle{empty}
\pagestyle{empty}

\begin{abstract}

The problem of finding minimizing  geodesics for a 
manifold $M$ with a sub-Riemannian structure is equivalent  
to the time optimal control of a driftless system on $M$ with a 
bound on the control. We consider here a class of 
sub-Riemannian problems on the classical Lie groups  $G$ 
where the dynamical equations are of the form 
$\dot x=\sum_j X_j(x) u_j$ and the $X_j=X_j(x)$ are right invariant 
vector fields on $G$ and $u_j:=u_j(t)$ the controls.   The vector fields $X_j$ are assumed to belong to the P part of a Cartan K-P decomposition. These types of problems admit a group of symmetries $K$ which act on $G$ by conjugation. Under the assumption that the minimal isotropy group \cite{Bredon}  in $K$ is discrete, we prove that we can reduce the problem  to a Riemannian problem on the regular part of the { associated  quotient space $G/K$}. On this part we define  the corresponding quotient metric. For the special cases of 
the K-P decomposition of $SU(n)$ of type {\bf AIII} we prove that the assumption  on the minimal isotropy group is verified. Moreover, under the assumption that the quotient space $G/K$ with the given metric has negative curvature we give a converse to a theorem of \cite{confran}, \cite{Dimitry}, proving that the cut locus has to belong to the singular part of $G$. As an example of applications of these  techniques we characterize the cut locus for a problem on $SU(2)$ of interest in the control of quantum systems.  

\end{abstract}

\vspace{0.5cm}


\vspace{0.5cm}

\section*{List of Symbols}
\noindent $T_x M$ -- Tangent space at $x$ for a manifold $M$\\
$TM$-- Tangent bundle of $M$, i.e., $TM:=\cup_{x \in M}T_xM$\\
$\Delta$ -- sub-bundle of the tangent bundle $TM$ \\
$\pi_{\Delta}:\Delta\to M$ -- restriction to $\Delta$ of the standard projection map from $TM$ to $M$ \\
$f_*|_x:T_xM\to T_{\pi(x)}N$ --  push-forward of a map $f:M\to N$ at the point $x\in M$ \\
$f_*:TM\to TN$ --  $f_*$ restricted to $T_xM$ is equal to $f_*|_x$\\
$G$ --  compact semisimple finite-dimensional real Lie group with corresponding Lie algebra $\mathfrak{g}$ identified with the tangent space at the identity. \\
$\bf{1}$ --  identity element of $G$ \\
$K=e^{\cal K}$ --  connected component containing ${\bf 1}$ of the Lie group associated to the Lie algebra $\cal K$ \\
$R_p:G\to G$ -- right multiplication  by $p$, $R_p(x)=xp$ \\
$L_p:G\to G$ -- left multiplication by $p$, $L_p(x)=px$ \\
$\text{ad}_P:\mathfrak{g}\to\mathfrak{g}$ -- adjoint map of a Lie algebra $\mathfrak{g}$ at $P\in\mathfrak{g}$, $\text{ad}_P(Q):=[P,Q]$ for every $Q\in\mathfrak{g}$ \\
$\langle\cdot|\cdot\rangle$ --  inner product induced by the Killing form\footnote{Or any positive scalar multiple of it.} of a Lie algebra $\mathfrak{g}$, $\langle P|Q\rangle:=-\text{tr}(\text{ad}_P\circ\text{ad}_Q)$ \\
$\langle \cdot,\cdot\rangle_{x}$ -- Riemannian metric at $x\in G$ or its sub-Riemannian restriction, $\langle R_{x*}P,R_{x*}Q\rangle_{x}:=\langle P|Q\rangle$ for $P,Q\in \mathfrak{g}$, \\
$\pi:G\to G/K$ -- quotient map associated to the action of the Lie group $K$ acting on the Lie group $G$ \\
$G_{\text{reg}}$ --  set of points in $G$ having minimal isotropy type.  \\
$G_{\text{sing}}$ -- set of non-regular points in $G$, i.e. $G_{\text{sing}}:=G - G_{\text{reg}}$ \\
$g_{\pi(x)}(\cdot,\cdot)$ --  Riemannian metric at $\pi(x)\in G_{\text{reg}}/K$.  \\
$d(p,q)$ --  sub-Riemannian distance from $p$ to $q$ in $G$. \\
$d_Q(\pi(p),\pi(q))$ --  Riemannian distance from $\pi(p)$ to $\pi(q)$ in $G_{\text{reg}}/K$ \\


\section{Introduction}  

Sub-Riemannian problems are equivalent to optimal control problems for driftless control systems when 
we want to minimize time with bounded energy or vice-versa  \cite{ABB}, \cite{confran}, \cite{Mont}. In these problems, one has a set of allowed directions at each point $p$ of a manifold $M$, with a given metric. One wants to transfer the state between two points  by moving  at each point following only the allowed directions and minimizing the corresponding distance (see next section  for formal definitions). In the paper \cite{Jurd}, V. Jurdjevi\'c introduced a class of  sub-Riemannian problems on matrix Lie groups $G$ 
 for which he was able to find an explicit expression of the optimal candidates. Such a class of problems, which were named K-P problems,  was then reconsidered in \cite{Bosca}, \cite{NoiAutomat}, \cite{confran}, because of their interest in quantum control. 
 In particular in \cite{NoiAutomat}, \cite{confran} an approach to their study 
 was used based on considering the symmetry action of a Lie subgroup of $G$, 
 $K \subseteq G$, on $G$.  This allowed the reduction of the number of unknown parameters in the optimal control law in several cases of interest. The action of $K$ on $G$ considered in \cite{NoiAutomat}, \cite{confran} is the {\it conjugation} (or {\it adjoint}) action where a matrix $x\in G$, is  transformed by a matrix $k \in K$ according to $x \rightarrow kxk^{-1}$. With this action the corresponding orbit space $G/K$ can often be mathematically described and visualized \cite{NoiAutomat} \cite{ADS}. However, in general,  such a space is not a manifold 
 but has the more general structure of a {\it stratified space} on which one has to generalize the standard notions of 
 differential geometry \cite{Sniat}. The strata are  (connected components of) the 
 orbit types  (see, e.g., \cite{Bredon}). Among them, the {\it minimal orbit type}, i.e., the orbit type corresponding to points with a minimal isotropy group, is, according  a theorem in the theory of Lie transformation groups,  an open and dense manifold in $G/K$, which is called the {\it regular part} { or {\it principal part}} of $G/K$ (or of $G$) \cite{Bredon}. The remaining part of $G/K$ (or of $G$) is called the {\it singular part}.

In this paper, we explore the possibility of studying K-P sub-Riemannian problems as Riemannian problems on the orbit space $G/K$.  Sub-Riemannian geodesics on $G$ can be obtained from Riemannian geodesics on $G/K$  as inverse images of the natural projection. We restrict ourselves to the regular part of $G/K$ and define a  
metric which allows us to obtain this reduction under the assumption that the minimum isotropy group in $K$ is discrete. 

Once the correspondence between a sub-Riemannian problem and a Riemannian one is established, one can use the powerful machinery of Riemannian geometry to answer questions for
 sub-Riemannian manifolds. We illustrate this in particular for the determination of the {\it sub-Riemannian cut locus}. The cut-locus is the set of points where geodesics lose optimality. It is of interest in Riemannian geometry for several reasons, including the fact that it gives information on the topology of the manifold (see, e.g., section 13.2 in \cite{DoCarmo}). In the optimal control context, the knowledge of the cut locus (from a given initial state $p$) is the first step to obtain the {\it complete optimal synthesis}, i.e, the knowledge  of all optimal trajectories. In fact, {\it all} the optimal trajectories are the ones with final point in the cut locus. Under the assumption that the quotient Riemannian manifold has negative curvature, we prove that the cut locus for the sub-Riemannian K-P problem has to be in the singular part of $G$, thus giving a converse to a theorem proved in \cite{confran} and \cite{Dimitry}. We present explicit calculations for an example in $SU(2)$.

The paper is organized as follows: In section \ref{Back}, 
 we give a more precise description of the K-P models 
 and their symmetries. In section \ref{metric} we describe in general 
 the choice of  the metric on the quotient space $G/K$ that 
 reduces the sub-Riemannian problem on $G$ to a Riemannian problem on $G/K$. Such a metric is well defined if the minimal isotropy group in $K$ is discrete.   In section \ref{AIII} we consider the K-P Cartan decomposition of $SU(n)$ of the type {\bf AIII} and prove that, in this case, this assumption is verified. The last two sections, 
 section \ref{CL} and \ref{CLSU2}, aim at showing how this approach can be used to determine the optimal synthesis for the sub-Riemannian problem and in particular to 
 compute the cut locus. In section \ref{CL} we establish the connection between the sub-Riemannian geodesics on $G$ and the Riemannian geodesics on the quotient space $G/K$ and investigate  the sub-Riemannian cut locus on $G$. We prove that if $G/K$ has negative curvature, the (Riemannian) cut locus on $G/K$, and therefore the (sub-Riemmanian) cut locus on $G$, has to belong to the singular part. In section \ref{CLSU2} we present explicit calculations (in coordinates) for an {\bf AIII} example on $SU(2)$. Such an example is of interest in quantum mechanics since it models the controlled evolution of a two level quantum system. We see that the curvature is negative in this case and therefore we characterize the cut locus using the above result.

\section{$K-P$ sub-Riemannian problems and their symmetries}\label{Back}

\subsection{Sub-Riemannian manifolds}

{ Given a Riemannian manifold $M$, a 
{\it sub-Riemannian structure} on $M$ is a 
subbundle $\Delta$ of the tangent bundle $TM$. Denoting by $\Delta_x$ the fiber at  
$x \in M$, (i.e., if $\pi_\Delta$ is the natural projection $\pi_\Delta \, : \, \Delta \rightarrow M$, $\Delta_x :=\pi_\Delta^{-1}(x)$) we assume that $\dim(\Delta_x)=m$ is constant, independent of $x$. In control theory $\Delta$ is often specified by giving a {\it distribution}, or {\it frame}, that is, a set of $m$ vector fields ${\cal F}:=\{ X_1,...,X_m\}$ such that, 
 for every $x \in  M$, $\texttt{span} \{X_1(x),...,X_m(x)\}= \Delta_x$. For simplicity we shall assume  that $\{ X_1,...,X_m\}$ is an orthonormal frame, that is, $\forall x \in M$, 
$\langle X_j(x), X_k (x) \rangle =\delta_{j,k}$, where $\langle \cdot , \cdot \rangle$ denotes an underlying Riemannian metric of $M$, and $\delta_{j,k}$ is the Kronecker delta. It is also assumed that the frame ${\cal F}$ is {\it bracket generating}, that is, if $Lie \, {\cal F}$ is the Lie algebra generated 
by ${\cal F}$, $Lie {\cal F}(x)=T_xM$, for every $x \in M$.   The $\Delta$ 
in this definition is frequently called a {\it horizontal distribution}.

A curve $\gamma:[a,b]\to M$ is said to be a {\it horizontal curve} if the tangent vector $\dot{\gamma}(t)\in \Delta_{\gamma(t)}$ for every $t\in[a,b]$.
For a manifold $M$, with a sub-Riemannian structure $\Delta$, one may 
define a metric structure:
\be{distB}
d(p,q):=\inf_{\gamma}\int_a^b\sqrt{\langle\dot\gamma(t), \dot\gamma(t)\rangle}dt
\ee
where $p,q\in M$ and the infimum is taken over all horizontal curves $\gamma:[a,b]\to M$ such that $\gamma(a)=p$, $\gamma(b)=q$. A horizontal curve $\gamma$ which is distance minimizing is called a {\it sub-Riemannian geodesic}. The above condition of ${\cal F}$ 
being  bracket generating 
guarantees, under the assumption that $M$ is a connected and complete metric 
space, that, for any two points $p$ and $q$, the distance is realized by a horizontal curve (cf. the Chow-Raschevskii theorem in \cite{ABB}, \cite{Mont}).  In the context of geometric control theory, sub-Riemannian manifolds arise in finite-dimensional smooth control systems with constraints. The horizontal distribution ${\cal F}$ describes the directions which may be taken at each point \cite{ABB}. The 
sub-Riemannian geodesics $\gamma$ parametrized by arclength 
($\|\dot \gamma (t) \|=c, \, \forall t $) are the trajectories which, for a given state transfer $p \rightarrow q$,  minimize time for the system 
$$
\dot x=\sum_{j=1}^mX_j(x)u_j(t), \qquad x(0)=p, 
$$ 
subject to the constraint $\|u\|^2 \leq c^2$ (cf., e.g., \cite{ABB}, \cite{confran}).

\subsection{$K-P$ problems}\label{KPPro}
 K-P sub-Riemannan problems were introduced in \cite{Jurd} and studied in \cite{NoiAutomat}, \cite{confran}, \cite{Bosca} in the context of quantum control. Such problems provide, in some sense, the simplest class of examples of sub-Riemannian 
problems (after the Riemannian case) since they are systems of {\it non-holonomy degree one} \cite{MurraySastry}: It is sufficient to take {\it one}  Lie bracket of the available vector fields to have a distribution which at every point spans the whole tangent space.  One main feature of these types of problems is that the equations of the Pontryagin Maximum principle of optimal control are integrable. Therefore an explicit expression of the optimal candidates can be obtained \cite{Bosca},  \cite{Jurd}. 

The basic setup of a K-P problem is as follows:
Let $G$ be a finite-dimensional, real, connected, semisimple Lie group with corresponding Lie algebra $\mathfrak{g}$. In the following, we will also 
take $G$ to  be compact, though many of the results below are generalizeable to 
the non-compact case. A {\it K-P Cartan decomposition} is a decomposition $\mathfrak{g}=\mathcal{K}\oplus\mathcal{P}$ into vector spaces $\cal K$ and $\cal P$ such that:
\be{KP}
[\mathcal{K},\mathcal{K}]\subseteq\mathcal{K} \qquad [\mathcal{K},\mathcal{P}] \subseteq \mathcal{P} \qquad  [\mathcal{P},\mathcal{P}] \subseteq \mathcal{K}. 
\ee
A bi-invariant, 
positive-definite, symmetric inner product $\langle\cdot|\cdot\rangle$ on $\mathfrak{g}$ (cf. \cite{Helg}, chapter III section 7, and chapter X) can be defined as follows: For $G$ compact, 
\be{Killing}
\langle P|Q\rangle:=-B(P,Q):=-\text{tr}(\text{ad}_P\circ\text{ad}_Q),
\ee 
the negative of the {\it Killing form}, {or any positive multiple of it}. The K-P problem is the minimum time problem  for  systems of the form 
\be{KPdrift}
\dot{X}=\sum_j u_j B_j X, \qquad X(0)={\bf 1}, 
\ee
with $X \in G$ and $u_j$ the controls, with $\|u \| \leq 1$. The elements $B_j$ form an orthonormal basis of ${\cal P}$ in a K-P decomposition and the sub-Riemannian structure is given by the frame of right invariant vector fields  
${\cal F}:=\{ B_1 X,..., B_m X\}$.  In other terms, the sub-bundle 
$\Delta$ is given by $\cup_{x \in G} R_{x*}{\cal P}$. Here we identify 
the tangent space of $G$ at the identity with the Lie algebra $\mathfrak{g}$, and let
 $R_x$ ($L_x$) denote the right (left) translation; The fiber at $x$ is given by $R_{x*} {\cal P}$ and it is spanned $\{ R_{x*}  B_j \}$. The Riemannian metric $\langle \cdot, \cdot \rangle$  on $G$ is derived from the above Killing form $B$ on 
 $\mathfrak{g}$, again 
 by identifying $\mathfrak{g}$ with the tangent space of $G$ at the identity: if $V_1$ and $V_2$ are tangent vectors at $x$ then $\langle V_1, V_2 \rangle:=-B(R_{x^{-1}*} V_1, R_{x^{-1}*} V_2)$. By definition the above  metric is right-invariant, i.e., $\forall x \in G$, $V_1, V_2$ in some tangent space of $G$, $\langle R_{x*} V_1,    R_{x*} V_2 \rangle =
 \langle  V_1,     V_2 \rangle$. By the properties of the Killing form, it also 
 follows that it is {\it left invariant} (same definition as before with $L_{x*}$ replacing $R_{x*}$) therefore the metric thus defined is {\it bi-invariant}.

Applying the Pontryagin maximum principle for the minimum time problem for system (\ref{KPdrift}) one finds,  \cite{Bosca}, \cite{Jurd},  
 that  the optimal 
 control has  the form $\sum_j B_ju_j(t)=e^{At}Pe^{-At}$ 
for some $A\in\cal K$ and $P\in\cal P$, and the corresponding optimal trajectory  has the form  $X(t)=e^{At}e^{(-A+P)t}$.

\subsection{Symmetries}

One of the main features of K-P problems is the existence of a {\it group of symmetries}. Consider the connected Lie group $K:=e^{\cal K}$ associated to the Lie algebra ${\cal K}$, and let us assume this Lie group to be compact.{\footnote{In the case with $G$ compact, this is guaranteed since $K$ is a closed Lie subgroup.}} This Lie group has a (left, proper) action  $\Phi_k$, on $G$ by conjugation, i.e, for $x \in G$, $k \in K$, 
 \be{Phikaction}
 \Phi_{k} x:=kxk^{-1}. 
 \ee
Such an action gives  a {\it symmetry} for the sub-Riemannian optimal control problem (\ref{KPdrift}) as it satisfies the following conditions: 1) For the initial condition,  which is the identity $\bf 1$ in (\ref{KPdrift}), and every $k$, $\Phi_k {\bf 1}={\bf 1}$, that is, the action leaves the initial condition unchanged 2) (invariance) 
For any $k \in K$, and given the distribution $\Delta$ defining the sub-Riemannian structure, we have $\Phi_{k*} \Delta_x=\Delta_{\Phi_k x}$. This property is a consequence of the fact that the $B_j$'s in (\ref{KPdrift}) form  a basis in ${\cal P}$ and from the second one in (\ref{KP}),  the set $\{ k B_j k^{-1} \, | \, j=1,2,...,m \}$ is a basis in ${\cal P}$ as well.\footnote{More in detail $\Phi_{k*} R_{x*} B_j:=L_{k*} R_{k^{-1}*} R_{x*} B_j=L_{k*} R_{k^{-1}*} R_{x*} R_{k*} R_{k^{-1}*} B_j=L_{k*}R_{kxk^{-1}*}R_{k^{-1}*}B_j=R_{kxk^{-1}*}L_{k*} R_{k^{-1}*}B_j$, since $L_{k*}$ and $R_{kxk^{-1}*}$ commute. However this is equal to $R_{\Phi_k(x)*} L_{k*} R_{k^{-1}*}B_j$ and since $L_{k*} R_{k^{-1}*}B_j\in {\cal P}$ it belongs to $\Delta_{\Phi_k(x)}$.} 3) For any $k \in K$, $\Phi_k$ is an {\it isometry}, that is, for any two tangent vectors $V$ and $W$ at $x \in G$, 
$$\langle \Phi_{k*}{V}, \Phi_{k*} W\rangle_{\Phi_kx}=\langle V, W\rangle_x,$$
where $\langle \cdot, \cdot \rangle_x$ is the Riemannian metric on $G$ 
calculated at $x$. In our case, the Riemannian metric is given by the Killing metric above described. We have $\langle \Phi_{k*} V, \Phi_{k*} W \rangle_{\Phi_k(x)}=\langle L_{k*} R_{k^{-1}*} V,      L_{k*} R_{k^{-1}*}  W \rangle_{\Phi_k(x)} =  \langle  V,    W \rangle_{x} $ because of the left and right invariance  (bi-invariance) of the metric.

A consequence of these properties is that (cf., \cite{confran})  if $\gamma(t)$ is a minimizing sub-Riemannian geodesic from the identity ${\bf 1}$ to $p \in G$, for any $k \in K$, $\Phi_k(\gamma(t))$  is a minimizing geodesic from ${\bf 1}$ to $\Phi_k(p)$. Therefore optimal geodesics are the `lifts' (cf. next section) of appropriate curves on the quotient space $G/K$ which, we will see, are also geodesics corresponding to an appropriate Riemannian metric.

The action of $K$ on $G$ by conjugation is proper (since $K$ is assumed to be compact) but it is not a {\it free} action since, for instance,  the isotropy group of the identity is the full group $K$. Therefore $G/K$ is not guaranteed to be a manifold and it is in fact a {\it stratified space}. To understand the stratified structure of    $G/K$, following the theory of Lie transformation groups (see, e.g., \cite{Bredon}), 
 one considers all the possible subgroups of $K$ which are isotropy groups for some elements in $G$. Groups $H$ which can be obtained one from the other by a similarity transformation ($H_2=kH_1 k^{-1}$), for $k \in K$ are placed in equivalence classes $(H)$ called {\it isotropy types}. Elements in $G$ which have isotropy groups with the same isotropy type $(H)$ are placed in the same set $G_{(H)}$,  and points on the same orbit must be in  the same orbit  type.  Therefore it makes sense to consider the quotient spaces $G_{(H)}/K$. The full quotient space $G/K$ is the disjoint union of the {\it orbit types} $G_{(H)}/K$'s  over all possible isotropy types $(H)$. The connected components of  $G_{(H)}/K$ are manifolds which are the strata in the stratified space $G/K$. Among the various isotropy types,  one can introduce a {\it partial ordering} by saying that $(H_1) \leq   (H_2)$ if and only if there exists a group in $(H_1)$ which is conjugate to a subgroup of 
a group in  $(H_2)$. The {\it minimum isotropy type theorem} (cf., e.g., \cite{Bredon}) states that there exists a minimum isotropy type $(H_{min})$ which is $\leq$ any isotropy type and the corresponding { orbit type} $G_{(H_{min})}/K$ is a {\it connected open and dense}  manifold in  $G/K$. It is called the 
{\it regular part} of $G/K$ and we shall denote it by $G_{\text{reg}}/K$. The remaining part $G/K-G_{\text{reg}}/K=:G_{\text{sing}}/K$ is called the {\it singular part}. Their preimage in $G$ are called the regular part $G_{\text{reg}}$ and singular part $G_{\text{sing}}$ of $G$, respectively. The following example which was also treated in \cite{NoiAutomat} clarifies this idea and will be the object of the analysis in section \ref{CLSU2}. 

\bex{PreesSU2}
Consider $G=SU(2)$ and the decomposition of $su(2)$ into diagonal matrices and antidiagonal matrices which give the ${\cal K}$ and ${\cal P}$ part of the Cartan K-P decomposition,  respectively. The Lie group $K:=e^{\cal K}$ is the (one dimensional) Lie group of diagonal matrices in $SU(2)$. Matrices that are diagonal in $SU(2)$ have as isotropy group the whole $K$ while matrices that are not diagonal have as isotropy group $\{ \pm {\bf 1} \}$. Therefore there exist only two isotropy types and the minimum isotropy type is given by the discrete group 
$\{\pm {\bf 1} \}$. Writing a general  element $X \in SU(2)$ as  
\be{generic}
X:=\begin{pmatrix} z & w \cr -w^* &z^*\end{pmatrix}, \qquad |z|^2+|w|^2=1,
\ee
conjugation by an element of    $K$ 
does not modify the diagonal element $z$, 
while it may arbitrarily change the phase of the antidiagonal element $w$. Therefore, 
 the orbits in $SU(2)/K$ are parametrized by the $(1,1)$ element, $z$, i.e., a point 
  in the closed unit disc of  the complex plane. If $|z|=1$ the isotropy group is $K$. This is the singular part of the orbit space. If $|z|<1$ then the isotropy group is $\{ \pm {\bf 1} \}$. This is the regular part corresponding to the interior of the unit disc. Section \ref{AIII} generalizes some of these properties. 
\eex

Symmetry reduction has a long history in control theory and we refer to \cite{dodici}, \cite{quattordici},  \cite{Martinez},  \cite{Tomizu}, as an entry point to an extensive literature. The novelty here is the application to K-P systems, the fact that the quotient space has a more general structure than the one of a manifold and the extensive use of Riemannian geometry.  

\section{A Riemannian metric on the quotient space} \label{metric}


We  define a Riemannian metric on the regular part of the quotient space, $G_{\text{reg}}/K$ as follows: Suppose the minimal isotropy type is discrete, {and recall that $G_{\text{reg}}$ is an open and dense submanifold of $G$ (c.f. \cite{Bredon}, Chapter IV Theorem 3.1)}. For $V,W\in T_{\pi(x)}(G_{\text{reg}}/K)$, let $P,Q\in\mathcal{P}$ such that $\pi_*R_{x*}P=V$ and $\pi_*R_{x*}Q=W$, where $\pi$ is the natural projection $\pi: \, G_{\text{reg}} \rightarrow G_{\text{reg}}/K$ 
and define the metric $g$ on $G_{\text{reg}}/K$, 
\be{metricdef}
g_{\pi(x)}(V,W):=\langle R_{x*} P, R_{x*} Q\rangle_x:=-\text{tr}(\text{ad}_P\circ\text{ad}_Q):=\langle P| Q \rangle, 
\ee
(cf., (\ref{Killing})). 
For this definition to be well posed, we must prove that the  `lifts' $P$ and $Q$ of $V$ and $W$, respectively,  exist and that the metric is independent of the choice of such lifts and the choice of basepoint $x$ in the fiber corresponding to $\pi(x)$. We first  observe that for $x\in G_{\text{reg}}$, $\pi_*:T_x G_{\text{reg}}\to T_{\pi(x)}(G_{\text{reg}}/K)$ is defined since $G_{\text{reg}}$ is an open dense submanifold of $G$ and therefore $T_x G_{\text{reg}}$ may be identified with $T_x G$. The projection $\pi$ in the definition and in the following is meant to be restricted to $G_{\text{reg}}$, so that $\pi_*$ is restricted to $T_xG_{\text{reg}}$.  

\bl{kerpi}
$\ker_x(\pi_*)=\{R_{x*}A-L_{x*}A|A\in\mathcal{K}\}$.
\el

\bpr
From the theory of proper actions of Lie transformation groups  
(e.g. \cite{DK}, chapter II), it follows that $\ker_x(\pi_*)=T_x(K \cdot x)$. Therefore, for a path $k(t):(-a,a)\to K$ with $k(0)={\bf 1}$, $X\in\ker_x(\pi_*)$ if and 
only if $X=\frac{d}{dt}|_{t=0}k(t)xk(t)^{-1}=R_{x*}\dot{k}(0)-L_{x*}\dot{k}(0)$.
\epr

\bl{PeqQ}
Suppose $P,Q\in\mathcal{P}$ and $x\in G_{\text{reg}}$ such that $\pi_*R_{x*}P=\pi_*R_{x*}Q$. Then $P=Q$.
\el

\bpr
$\pi_*R_{x*}(P-Q)=0$ if and only if $R_{x*}(P-Q)\in\ker_x\pi_*$, and therefore by Lemma \ref{kerpi}, there exists an $A\in\mathcal{K}$ such that $R_{x*}(P-Q)=R_{x*}A-L_{x*}A$. Applying $R_{x^{-1}*}$ to both sides yields\footnote{Here we adopt the common abuse of notation that $L_{x*}R_{x^{-1}*}A=xAx^{-1}$.}: $P-Q=A-xAx^{-1}$. Since  $P-Q\in\mathcal{P}$ and $A\in\mathcal{K}$ we have  that the $\mathcal{K}$-part of $xAx^{-1}$ is $(xAx^{-1})_{\mathcal{K}}=A$, while the $\mathcal{P}$-part is $(xAx^{-1})_{\mathcal{P}}=Q-P$. Therefore, since $\mathcal{K}$ and $\mathcal{P}$ are orthogonal with respect to the Killing metric (\ref{Killing}), we have that:
\be{AxA}
 \langle A|xAx^{-1}\rangle=\langle A|A\rangle
\ee
and
\be{QxA}
\langle Q|xAx^{-1}\rangle=\langle Q|Q-P\rangle
\ee
Now, we will compute $\langle P|P\rangle$, then rearrange to show that $\langle P-Q|P-Q\rangle=0$, which because of positive semidefiniteness of $\langle\cdot|\cdot\rangle$ will imply that $P-Q=0$.
\be{PP}
\begin{split}
\langle P|P\rangle=\langle Q+A-xAx^{-1}|Q+A-xAx^{-1}\rangle= \\
=\langle Q|Q\rangle+2\langle Q|A\rangle-2\langle Q|xAx^{-1}\rangle+\langle A|A\rangle \\
-2\langle A|xAx^{-1}\rangle+\langle xAx^{-1}|xAx^{-1}\rangle
\end{split}
\ee
Now, applying the orthogonality of $\mathcal{P}$ and $\mathcal{K}$ to $\langle Q|A\rangle$; applying (\ref{QxA}) to $\langle Q|xAx^{-1}\rangle$; applying (\ref{AxA}) to $\langle A| xAx^{-1}\rangle$; and applying the bi-invariance of the metric to $\langle xAx^{-1}|xAx^{-1}\rangle=\langle A|A\rangle$ yields:
\be{PP2}
\begin{split}
\langle P|P\rangle=\langle Q|Q\rangle-2\langle Q|Q-P\rangle \\
=2\langle Q|P\rangle-\langle Q|Q\rangle
\end{split}
\ee
Rearranging yields $\langle P|P\rangle-2\langle P|Q\rangle+\langle Q|Q\rangle=0$, and so $\langle P-Q|P-Q\rangle=0$.
\epr

\bl{PpQq}
Suppose $p,q\in G_{\text{reg}}$, $k \in K$ such that $p=kqk^{-1}$ and $P,Q \in \mathcal{P}$ with $\pi_*R_{p*}P=\pi_*R_{q*}Q$. Then $P=kQk^{-1}$.
\el
\bpr
\be{RqQ}
\pi_*R_{q*}Q=\pi_*R_{p*}P=\pi_*R_{kqk^{-1}*}P=\pi_*R_{k^{-1}*}R_{q*}R_{k*}P
\ee
where, in this last line, we have applied the antihomomorphism property of right multiplication. Now, since $\pi$ is constant on equivalence classes, $\pi_*\circ(L_{k^{-1}*}R_{k*})=\pi_*$, and so (\ref{RqQ}) becomes:
\be{RqQ2}
\pi_*R_{q*}Q=\pi_*L_{k^{-1}*}R_{k*}R_{k^{-1}*}R_{q*}R_{k*}P
\ee
$$
=\pi_*L_{k^{-1}*}R_{q*}R_{k*}P
$$
Then since the left multiplication operator and right multiplication operators commute this yields:
\be{RqQ3}
\pi_*R_{q*}Q=\pi_*R_{q*}(k^{-1}Pk)
\ee
So by Lemma \ref{PeqQ}, $Q=k^{-1}Pk$.
\epr

Taken together, these lemmas imply that if, given $V, W\in T_{\pi(x)}G_{\text{reg}}/K$, one is able to find $P,Q\in\mathcal{P}$ such that $\pi_*R_{x*}P=V$ and $\pi_*R_{x*}Q=W$, then (\ref{metricdef}) is well-defined. Lemma \ref{PeqQ} implies that if such $P,Q$ exist, then once the basepoint $x$ in $\pi^{-1}(\pi(x))$ has been fixed, $P,Q$ are unique. Lemma \ref{PpQq} implies that if one lifts using a different basepoint, say $y=kxk^{-1}$ for some $k\in K$, then the corresponding $P,Q$ will be $kPk^{-1},kQk^{-1}$, and so by the bi-invariance of the metric, $\langle kPk^{-1}|kQk^{-1}\rangle=\langle P|Q\rangle$, and so the definition (\ref{metricdef}) will not depend on the choice of the basepoint. The question then becomes: When does such a lift exist?
\bt{Metric}
The  metric  in (\ref{metricdef}) is defined  if and only if 
the minimal isotropy type in $K$ is discrete.
\et
\bpr
Fix $x\in G_{\text{reg}}$ and consider the linear map 
\be{linmap}
\pi_*R_{x*}:\mathcal{P}\to T_{\pi(x)}G_{\text{reg}}/K
\ee
By the above remarks, it suffices to prove that $\pi_*R_{x*}$ is an isomorphism. Since Lemma \ref{PeqQ} implies that (\ref{linmap}) is injective, the theorem is 
proved if we show that the map is surjective. This is the case if and only if $\dim\mathcal{P}=\dim(T_{\pi(x)}(G_{\text{reg}}/K))=\dim(G_{\text{reg}}/{K})$. On the other hand, $\dim(G_{\text{reg}}/K)=\dim G-\dim K+\dim H_{\text{min}}$, where $H_{\text{min}}$ is any subgroup of $K$ which is of minimal isotropy type (see, e.g., Theorem 2.3 in \cite{Dimitry}). Since $\dim G=\dim\mathcal{P}+\dim\mathcal{K}$, we have that $\pi_*R_{x*}$ is surjective if and only if $\dim H_{\text{min}}=0$, i.e. if and only if $H_{\text{min}}$ is discrete.
\epr
\br{Rema1}
An alternative to the definition (\ref{metricdef}) could have been  to see the projection $\pi: G_{\text{reg}} \rightarrow G_{\text{reg}}/K$ as a Riemannian submersion \cite{DoCarmo} and define a {\it vertical} distribution given by $\ker_x\pi_*$ (described in Lemma \ref{kerpi}) at any point $x$ and the {\it horizontal} distribution given by the orthogonal space (to the vertical one) in the given Riemannian metric. A metric is defined analogously to what we have done here taking for each tangent vector in the quotient space its `lift' to the horizontal space. This however does not coincide with the ${\cal P}$ space of the sub-Riemannian structure.   
\er 
An important feature of the above defined metric is that the length is preserved going from horizontal curves $\gamma$ in $G_{\text{reg}}$ to the corresponding curves $\pi(\gamma)$ in $G_{\text{reg}}/K$. If $\gamma$ is a horizontal curve in $G_{\text{reg}}$, then from (\ref{metricdef}) and since $\dot \gamma(t) \in R_{\gamma(t)*}{\cal P}$ we have $g_{\pi(\gamma(t))}(\pi_* \dot \gamma(t), \pi_* \dot \gamma(t))=
\langle \dot \gamma(t), \dot \gamma(t) \rangle_{\gamma(t)}$ and therefore from (\ref{distB}) the length is preserved. If $\gamma=\gamma[0,T]$ has one of the endpoints in $G_{\text{sing}}$, the length is preserved on any sub-interval of $[0,T]$. This is the case of interest for us since our initial point, the identity ${\bf 1}$, belongs to the singular part of $G$.


\section{$K-P$ Problems of the type {\bf AIII} on $SU(n)$} \label{AIII}
A K-P Cartan decomposition of $su(n)$ of the type {\bf AIII} is a 
decomposition $su(n)={\cal K} \oplus {\cal P}$ satisfying (\ref{KP}) where ${\cal K}$ consists of block diagonal matrices in $su(n)$ and ${\cal P}$ are block anti-diagonal matrices. More specifically let $q \leq n-q$.  Then the matrices in ${\cal K}$ have the form 
$\begin{pmatrix} A_{q \times q} & {\bf 0} \cr {\bf 0} & B_{(n-q) \times (n-q)}\end{pmatrix}$ with 
$Tr(A_{q \times q})+Tr(B_{(n-q) \times (n-q)})=0$,  while the matrices in ${\cal P}$ 
have the form  $\begin{pmatrix} {\bf 0}_{q \times q} & C_{q \times {(n-q)}} \cr -C_{q \times {(n-q)}}^\dagger  & {\bf 0}_{(n-q) \times (n-q)}\end{pmatrix}$. The Lie group $K$ is the Lie subgroup of $SU(n)$ of block diagonal matrices with blocks of dimension $q \times q$ and $(n-q) \times (n-q)$. If we consider the left (or right) multiplication action of $K$ on $SU(n)$, the quotient space is one of the symmetric spaces classified by Cartan \cite{Helg}. If we consider the conjugation action as we do in this paper, the quotient space is one of the stratified spaces discussed above. We show here 
that in this case we can define a metric as in theorem \ref{Metric} of the   above section. 

\begin{theorem} Consider a K-P Cartan decomposition of the type {\bf AIII}.  
The minimum isotropy group in $K$ for the conjugation action on $SU(n)$ is the discrete (Abelian) group $H:=\{{\bf 1}, \omega  {\bf 1}, ..., \omega^{n-1}  {\bf 1} \}$, where 
$\omega:=e^{i\frac{2 \pi}{n}}$. 
\end{theorem}
In particular, the K-P problem of the type 
${\bf AIII}$ satisfies the condition of Theorem \ref{Metric} to define 
a quotient metric on the regular part.

\bpr
It is enough to show that there is a matrix in $SU(n)$ whose 
isotropy group in $K$ is $H$. In order to do that, 
we set up  a few definitions. 

Let $n:=kq+j$, with $k \geq 2$ and $j=0,1,...,q-1$, and consider matrices in $SU(n)$ represented in blocks: $k$ block-rows with $k$ $q \times q$ blocks followed 
by a $q \times j$ block plus a last $k+1$-th block-row with $k$ $j \times q$ blocks 
followed by  a $j\times j$ block. In this setting, we define,  for $l=1,...,k-1$, the matrices $A_l$ which are equal to the identity except in the, $q \times q$,  
blocks $(l,l)$, $(l,l+1)$, $(l+1,l)$, $(l+1, l+1)$,  which are occupied by 
$\begin{pmatrix} \frac{1}{\sqrt{2}} {\bf 1}_{q} &  
\frac{1}{\sqrt{2}} {\bf 1}_{q} \cr  - \frac{1}{\sqrt{2}} {\bf 1}_{q} & \frac{1}{\sqrt{2}} {\bf 1}_{q} \end{pmatrix}$. We also define for $l=1,...,k-2$ the matrices 
$\hat P_l$ (if any) which have a special unitary $q \times q$ matrix $P_l$ in the $(l,l)$ block and are elsewhere  equal to the identity. The matrix $\hat P_{k-1}$ is the identity except for the last three block rows and columns, i.e., the blocks $B_{r,s}$, $r,s=k-1,k,k+1$, which are occupied by a special unitary matrix 
$\begin{pmatrix} P_{k-1} &  0 & 0 \cr 0 & T_{1,1} & T_{1,2} \cr 
0 & T_{2,1} & T_{2,2} \end{pmatrix}$,  
of dimensions $(2q+j) \times (2q+j)$. The matrices $T_{1,2}$, $T_{2,1}$ and $T_{2,2}$ only exist if $j\geq 1$. Here we choose the matrices as follows: 

\begin{enumerate}
\item (R1) $P_{k-1}$ a diagonal matrix in $U(q)$ with all the elements on the diagonal different from each other. 

\item (R2) $T_{1,1}$ is a $q\times q$ matrix whose first row is such that all non diagonal elements (if any) are $\not=0$. 

\item (R3) $T_{1,2}$ (if existing), is a $q \times j$ matrix having a generalized left inverse, i.e., a $j \times q$ matrix $L$ such that $LT_{1,2}=\lambda {\bf 1}_j$, for a scalar $\lambda \not=0$.   

\end{enumerate}

 For $l=1,..,k-2$ define (if any) the matrices 
\be{effel}
F_l:=A_l \hat P_l, 
\ee
and for $l=k-1$ define 
\be{effeelle1}
F_{k-1}:=A_{k-1} \hat P_{k-1}. 
\ee
The matrix we will show to have isotropy group given by the scalar matrices in $K$ will have the form $F_1F_2\cdots F_{k-2} \hat F_{k-1}$ for an appropriate choice of the above matrices $P_1,...,P_{k-2}$, and with the requirements (R1)-(R3). Let $R$ be an element in the isotropy group of  $F_1F_2\cdots F_{k-2} \hat F_{k-1}$ in $K$. From the block diagonal form of $R$, we write it as 
\be{erre}
R:=\begin{pmatrix} Q & 0 & \cdots & 0 \cr 
0 & Y_{1,1} & \cdots & Y_{1,k} \cr 
0 & Y_{2,1} & \cdots & Y_{2,k} \cr 
\cdot & \cdot & \cdot & \cdot \cr 
\cdot & \cdot & \cdot & \cdot \cr
\cdot & \cdot & \cdot & \cdot \cr
0  & Y_{k,1} & \cdots & Y_{k,k}, \end{pmatrix}
\ee
 for a unitary matrix $Q$ of dimension $q \times q$. If $R$ is in the isotropy group of $F_1F_2\cdots F_{k-2} \hat F_{k-1}$, we must have 
\be{tocomm}
RF_1F_2\cdots F_{k-2} \hat F_{k-1}=F_1F_2\cdots F_{k-2} \hat F_{k-1}R. 
\ee
Using (\ref{erre}) and the definition of $F_1,...,F_{k-2}, \hat F_{k-1}$, by comparing the first column-block  on the left and right hand side of (\ref{tocomm}) we get,
\be{tobus3}
QP_1=P_1 Q, \qquad Y_{1,1} P_1=P_1 Q, 
\ee
which implies in (\ref{erre}) that $Y_{1,1}=Q$ and therefore $Y_{1,2}$,...,$Y_{1,k}$ are all zero as well as $Y_{2,1}$ through $Y_{k,1}$. Therefore $R$ must have the form 
\be{errenew}
R:=\begin{pmatrix} Q & 0 & 0& \cdots & 0  \cr 
0 & Q &  0& \cdots & 0  \cr 
0 & 0 & Y_{2,2} &  \cdots & Y_{2,k}  \cr 
\cdot & \cdot & \cdot & \cdot & \cdot  \cr 
\cdot & \cdot & \cdot & \cdot & \cdot  \cr
\cdot & \cdot & \cdot & \cdot & \cdot \cr
0  & 0 & Y_{k,2} & \cdots & Y_{k,k}. \end{pmatrix}. 
\ee 
 The matrix $R$ with this form commutes with $F_1$ because of (\ref{tobus3}) 
 and therefore (\ref{tocomm}) reduces to 
\be{tocomm2}
RF_2\cdots F_{k-2} \hat F_{k-1}=F_2\cdots F_{k-2} \hat F_{k-1}R. 
\ee
Comparing the second block columns of the left and right hand side we get, similarly to 
(\ref{tobus3}) 
\be{tobus4}
QP_2=P_2 Q, \qquad Y_{2,2} P_2=P_2 Q, 
\ee
which implies $Y_{2,2}=Q$ in (\ref{errenew}), with zero in the remaining places in the corresponding block-row and block-column in $R$, and moreover $R$ commutes with $F_2$. Continuing this way until $k-2$, we obtain that $R$ has the block diagonal form 
$R=\text{diag} (Q,Q,...,Q, \hat R)$ for $k-1$  blocks $Q$ (of dimensions $q \times q$ 
and a block $\hat R$ of dimension $(q+j) \times (q+j)$.  Moreover formulas (\ref{tocomm}),  (\ref{tocomm2}) reduce to $R  \hat F_{k-1}=\hat F_{k-1} R$, which for the last $(2q+j) \times (2q+j)$ rows and columns gives 
\be{lastrc}
\begin{pmatrix} 
Q & 0 & 0 \cr 
0 & Y_{k-1,k-1} & Y_{k-1,k} \cr 
0 & Y_{k, k-1} & Y_{k,k} 
\end{pmatrix}
\begin{pmatrix}
P_{k-1}  & T_{1,1} & T_{1,2} \cr 
-  P_{k-1} & T_{1,1} & T_{1,2} \cr 
0 & T_{2,1} & T_{2,2}
\end{pmatrix} 
=
\ee
$$
\begin{pmatrix}
P_{k-1}  & T_{1,1} & T_{1,2} \cr 
-  P_{k-1} & T_{1,1} & T_{1,2} \cr 
0 & T_{2,1} & T_{2,2}
\end{pmatrix} 
\begin{pmatrix} 
Q & 0 & 0 \cr 
0 & Y_{k-1,k-1} & Y_{k-1,k} \cr 
0 & Y_{k,k-1} & Y_{k,k} 
\end{pmatrix}. 
$$
From equality of the first block column we get $Y_{k-1,k-1}=Q$. The proof is already finished in the case $q=1$. The matrix $R$ is a scalar matrix in $SU(n)$ in this case. If $q\not=1$ and $j\not=0$,   $Y_{k-1,k-1}=Q$ implies that $Y_{k,k-1}$ and $Y_{k-1,k}$ are zero. Moreover,  independently of the value of $j$, the requirement (R1) for $P_{k-1}$ above and equality of the first block-column give that 
$Q=Y_{k-1,k-1}$ is a diagonal matrix. Consider now equality of the second block column using  the fact that   $Q=Y_{k-1,k-1}$ is diagonal. This gives  that $T_{1,1}$ commutes with $Q$ and since $Q$ is diagonal and for the requirement (R2) above on $T_{1,1}$,  $Q$ is a scalar matrix $Q:=e^{i\phi} {\bf 1}_q$. The proof is finished if $j=0$. Otherwise, consider the third block-column and in particular equality of the first block. It gives $QT_{1,2}=T_{1,2} Y_{k,k}$. Using the requirement (R3) above  and multiplying by $L$ and using the fact that $Q=e^{i\phi} {\bf 1}_q$, we get $Y_{k,k}=e^{i\phi} {\bf 1}_j$. Therefore $R$ in (\ref{erre}) must be a scalar matrix. Since it is in $SU(n)$ it must belong to the discrete (in fact finite) group $H$. This  completes the proof. 
 \epr  

{
\section{The sub-Riemannian and Riemannian geodesics and the cut locus}\label{CL}

From now on, we will assume that the metric described in Section \ref{metric} exists.
The form of the optimal sub-Riemannian geodesics is known \cite{Bosca}, \cite{Jurd},  and it was described at the end of subsection \ref{KPPro}. In particular,  optimal geodesics are {\it analytic} and one can apply Corollary 3.6 in \cite{confran}\footnote{Notice that there is a mismatch of terminology between \cite{confran} and this paper. The 
cut locus in \cite{confran} is called `{\it critical locus}' while the name `cut locus' is used for points where  two or more minimizing geodesics intersect.}: If a point 
 $q$ is  in $G_{\text{reg}}$, then the optimal sub-Riemannian geodesic connecting  a point $p$ to  $q$ is entirely contained in   $G_{\text{reg}}$, except, possibly, for the initial point $p$. Therefore when optimally connecting two points in $G_{\text{reg}}$,  we can restrict ourselves to curves entirely contained in $G_{\text{reg}}$. In the following, we shall denote by $d(\cdot, \cdot)$ the sub-Riemannian distance on $G$ and therefore $G_{\text{reg}}$ defined by (\ref{distB}) and by 
$d_Q(\cdot, \cdot)$ the Riemannian distance 
on $G_{\text{reg}}/K$ with the metric 
defined in section \ref{metric}. If two points $p$ and $q$ are in $G_{\text{reg}}$, since the length is preserved by the projection $\pi$ under the adopted metric, as discussed at the end of section \ref{metric},  we have 
\be{inemetric}
d_{Q}(\pi(p), \pi(q)) \leq d(p,q). 
\ee 
The following theorem establishes the connection between Riemannian geodesics in $G_{\text{reg}}/K$ and sub-Riemannian geodesics in $G$, starting from the identity ${\bf 1}$.  
\bt{BasicTheo}
Assume $\gamma=\gamma(t)$ is a sub-Riemannian geodesic defined in $[0,T]$ optimally connecting ${\bf 1}$ and $q \in G_{\text{reg}}$. Then $\pi(\gamma)$ is a Riemannian geodesic from $\pi(\gamma(t_0))$ to 
$\pi(\gamma(T))=\pi(q)$, for any $t_0 \in (0,T)$. 

Moreover 

\be{state2}
\lim_{t_0 \rightarrow 0^{+}} d_{Q}(\pi(\gamma(t_0)), \pi(q))=d({\bf 1}, q). 
\ee
\et 

\bpr
Assume that there exists a $t_0 \in (0,T)$ such that the geodesic between $\pi(\gamma(t_0))$ and $\pi(q)$ is not $\pi(\gamma)$. Denote such a geodesic by $\Gamma$. 
Let $\bar t $ the smallest value of $t \geq t_0$  such that $\Gamma(t)=\pi(\gamma(t))$ for $t=\bar t$ and 
\be{different}
\Gamma(t)\not=\pi(\gamma(t)), \forall t \in (\bar t, \bar t+ \bar \epsilon), 
\ee 
for some appropriate $\bar \epsilon > 0$. 

Recall that from the proof of Theorem \ref{Metric} that for any $x \in G$, $\pi_{*}|_x$ { restricted to $R_{x*}\mathcal{P}$} is an isomorphism from $R_{x*} {\cal P}$ to $T_{\pi(x)}G_{\text{reg}}/K$, denoting by $\pi_{*}|_x^{-1}$ its inverse, consider (in local coordinates) the differential equation 
\be{lifteq}
\dot \gamma_1=\pi_{*}|_{\gamma_1(t)}^{-1} \dot \Gamma(t),  \qquad \gamma_1(\bar t)=\gamma(\bar t), 
\ee
which has a unique solution $\gamma_1$ in $[\bar t -\epsilon, \bar t +\epsilon]$ for appropriate $\epsilon$, choosing $\epsilon < \bar{\epsilon}$. Moreover $\pi(\gamma_1)=\Gamma$. 

Denote by $L$ the length of $\gamma_1$ between $\gamma_1(\bar t-\epsilon)$ and $\gamma_1(\bar t+\epsilon)$, which is $\geq $ the (sub-Riemannian) distance 
$d(\gamma_1(\bar t-\epsilon),\gamma_1(\bar t+\epsilon))$. Since $\pi$ preserves the distance we have 
$L=d_Q(\pi(\gamma_1(\bar t-\epsilon), \pi (\gamma_1(\bar t + \epsilon)) \leq d(\gamma_1(\bar t- \epsilon),\gamma_1(\bar t+ \epsilon))$  because of (\ref{inemetric}). Therefore $\gamma_1$ is a geodesic between $\gamma_1(\bar t-\epsilon)$ and $\gamma_1(\bar t+\epsilon)$. Moreover $\gamma(t)$ coincides with $\gamma_1(t)$, for $t \in [\bar t-\epsilon , \bar t]$. Since they are both geodesics and coincide on an open interval    $(\bar t-\epsilon , \bar t)$, because of analyticity of geodesics, they must coincide, which contradicts   (\ref{different}).

The above proof also shows that for every $t_0$ 
$$
d(\gamma(t_0), q)=d_Q(\pi(\gamma(t_0)), \pi(q))
$$
Taking the limit when $t_0 \rightarrow 0$ and using the continuity of the distance function $d$ from the Chow-Rashevski theorem we obtain (\ref{state2}). 
\epr
}
The theorem suggests a way to calculate the sub-Riemannian geodesics to points $q$ in $G_{\text{reg}}$ using Riemannian geometry. One calculates Riemannian geodesics $\Gamma$ leading to $\pi(q)$ in $G_{\text{reg}}/K$ and then calculate the `lift' i.e. the sub-Riemannian geodesic $\gamma_1$ leading to $q$ such that $\pi(\gamma_1)=\Gamma$ (cf. (\ref{lifteq})). Our main use of this correspondence is in the determination of the sub-Riemannian cut locus in $G$. In fact, from the above recalled results of \cite{confran}, if a geodesic intersects the singular part $G_{\text{sing}}$ of $G$, crossing the regular part,  then the point of intersection belongs to the cut locus. Taking into account 
that  $G_{\text{reg}}$ is open and dense in $G$, this suggests that the points in $G_{\text{sing}}$ are good candidates to belong to the sub-Riemannian cut locus. We prove next that, under the assumption of  {\it negative sectional curvature} of $G_{\text{reg}}/K$ the whole cut locus belongs to $G_{\text{sing}}$.

\bt{negativo}
Suppose that the sectional curvature of $G_{\text{reg}}/K$ under the metric of Section \ref{metric} is nonpositive and that $G_{\text{reg}}/K$ is simply connected. Then the intersection of the cut locus of $\{\bf{1}\}$ with $G_{\text{reg}}$ is empty.
\et
\bpr
{Assume $q$ is a cut point in $G_{\text{reg}}$ and $\gamma$ the corresponding sub-Riemannian geodesic (parametrized by constant speed) defined in $[0,T]$, with $\gamma(0)={\bf 1}$ and $\gamma(T)=q$. Then according to 
Theorem \ref{BasicTheo} $\pi(\gamma(t))$ is a minimimizing geodesic from $\pi(\gamma(t_0))$ to $\pi(q)$, for every $t_0 \in (0, T)$. Let $t_1>T$ sufficiently small so that $p:=\gamma(t_1)$ is still in $G_{\text{reg}}$ and the extension of  $\pi(\gamma)$ to $(0,t_1]$ is still a {\it locally minimizing} geodesic in $G_{\text{reg}}/K$ (with constant speed).\footnote{It coincides with $\pi(\gamma)$ in $(0,T]$ and therefore it satisfies the same geodesic equations \cite{DoCarmo} and initial conditions 
at $T$.}  Since $q=\gamma(T)$ is a cut point there exists   another sub-Riemannian optimal geodesic $\eta$ joining ${\bf 1}$ with $p=\gamma(t_1)$, and clearly $\pi(\eta)$ is a (locally) minimizing geodesic in $(0,t_1]$. Now, the theorem is proved if we prove the following claim:   

}

{\bf Claim } {\it Two locally minimizing geodesics $\hat \gamma$ and $\hat \eta$ in $G_{\text{reg}}/K$ such that  
$\lim_{t\rightarrow  0^+} \hat \gamma(t)= \lim_{t\rightarrow 0^+} \hat \eta(t) = \pi(\bf{1})$ cannot intersect in $G_{\text{reg}}/K$.}

Suppose 
$\hat \gamma(t_1)=\hat \eta(t_1)=\hat p\in G_{\text{reg}}/K$ (one may always reparametrize one of the geodesics so that they intersect at the same time $t_1$). Define the continuous, nonnegative,  function 
$f:[0,t_1]\to\mathbb{R}$ by $f(t):=d_Q(\hat \gamma(t), \hat \eta(t))$ for $t \in (0,t_1]$ and $f(0):=0$; 
Since $G_{\text{reg}}/K$ is a simply connected, complete smooth Riemannian manifold with nonpositive sectional curvature, it  is a {\it Hadamard manifold}.\footnote{We take this to be the definition of a Hadamard manifold. Other equivalent definitions exist (see  \cite{Ballmann}, Proposition 5.1).}. Therefore, $f(t)$ is a convex function (see  \cite{Ballmann}, Proposition 5.4). Therefore, for every $t\in [0,1]$, we have:
\be{conv}
0\leq f(t\cdot t_1)=f(t\cdot t_1+(1-t)\cdot 0)\leq tf(t_1)+(1-t)f(0)=0
\ee
So, $f(t\cdot t_0)=0$ for every $t\in [0,1]$, implying $\hat \gamma(t) = \hat \eta(t)$ for 
every $t\in (0,t_1]$. Therefore, two different geodesics starting from $\pi(\bf{1})$ cannot intersect in 
$G_{\text{reg}}/K$, and the Claim and therefore the theorem  are proved. 
\epr 

}

\section{Example: The cut locus for $SU(2)$}\label{CLSU2} 
K-P problems on $SU(n)$ are particularly interesting because $SU(n)$ may represent quantum mechanical evolutions of $n-$level quantum systems. In this context, the time optimal control problem is especially motivated because of the need to obtain fast computations in quantum information and because fast evolution is a way to avoid the degrading of the 
quantum evolution due to the effect of the environment, the so-called {\it decoherence}. Furthermore geometric time optimal control theory is a way to study the fundamental limitations of quantum evolution, the so-called {\it quantum speed limit}, and the related time-energy uncertainty relations (see, e.g., \cite{Brasil} and the references therein).  The case of $SU(2)$ is the simplest one but also a very important one as it models the evolution of two-level quantum systems, quantum bits, which are the basic building block of quantum information in the circuit based model \cite{NC}. This example has been treated in detail in \cite{NoiAutomat} where a method to find the time 
optimal control law from the identity to any final condition was described. The description of the cut locus was somehow implicit in \cite{NoiAutomat}. We derive it here with the methods of this paper.

For $G=SU(2)$ and $K=e^{\cal K}$ the (Abelian) Lie subgroup  of diagonal matrices in $SU(2)$, we have seen in Example \ref{PreesSU2} that $G/K$ is homeomorphic to the closed unit disc in the complex plane, and that $G_{\text{reg}}/K$ is the open unit disc. Explicitly, the homeomorphism is given by mapping a matrix in $SU(2)$ to its $(1,1)$-entry $z$ as in (\ref{generic}). 
By setting $z=x+iy$, we will use $x,y$ as coordinates in the open unit disc $G_{\text{reg}}/K$. Now, in order to compute the components $g_{ij}$ of the metric (\ref{metricdef}) in these coordinates, we need to know, at a point $z$ in the open disc, what $\frac{\partial}{\partial x}$ and $\frac{\partial}{\partial y}$ may be lifted to in the fiber $\pi^{-1}(z)$. So, lifting to a point (cf. (\ref{generic})) 
$q=\begin{pmatrix}
z & w \cr 
-w^* & z^*
\end{pmatrix} \in SU(2)$,
we would like to find a matrix $P=\begin{pmatrix}
0 & a+bi \cr
-a+bi & 0
\end{pmatrix}\in\mathcal{P}$ such that:
\be{SU2lift}
\pi_*R_{q*}P=\pi_*\left(
\begin{pmatrix}
0 & a+bi \cr
-a+bi & 0
\end{pmatrix}
\begin{pmatrix}
z & w \cr
-w^* & z^*
\end{pmatrix}
\right)
=\frac{\partial}{\partial x}
\ee
Letting $w=w_R+iw_I$, this implies:
\be{SU2lift2}
\begin{pmatrix}
-aw_R-bw_I \cr
aw_I-bw_R
\end{pmatrix}
=
\begin{pmatrix}
1 \cr
0
\end{pmatrix}
\ee
Therefore, $a=\frac{-w_R}{1-|z|^2}$, $b=\frac{-w_I}{1-|z|^2}$. Similarly, one may find a $Q=\begin{pmatrix}
0 & c+di \cr
-c+di & 0
\end{pmatrix}\in\mathcal{P}$ so that $\pi_*R_{q*}Q=\frac{\partial}{\partial y}$; in this case $c=\frac{w_I}{1-|z|^2}$, $d=\frac{-w_R}{1-|z|^2}$. Using this and the definition (\ref{metricdef}) 
(with the Killing metric $-\text{tr}(\text{ad}_P\circ\text{ad}_Q)= -\frac{1}{2}\text{Tr}(AB)$)
we have that the components of the metric on the regular part of the quotient space are given by $g_{ij}(z)=\frac{1}{1-|z|^2}\delta_{ij}$ with $i,j\in\{x,y\}$. Recalling that $z=x+iy$ and letting $r^2=|z|^2=x^2+y^2$, we may compute the Christoffel symbols $\Gamma^i_{kl}= \frac{1}{2}\sum_m g^{im}(\frac{\partial g_{mk}}{\partial x^l}+\frac{\partial g_{ml}}{\partial x^k}-\frac{\partial g_{kl}}{\partial x^m})$, using the standard formulas (in the case of a Riemannian connection; cf., e.g., \cite{DoCarmo} formula (10) Chapter 2, section 3)  at the point $(x,y)$:
\be{chris}
\begin{pmatrix}
\Gamma^x_{xx}=\frac{x}{1-r^2}, & \Gamma^x_{xy}=\Gamma^x_{yx}=\frac{y}{1-r^2}, & \Gamma^x_{yy}=\frac{-x}{1-r^2} \cr
\Gamma^y_{xx}=\frac{-y}{1-r^2}, & \Gamma^y_{xy}=\Gamma^y_{yx}=\frac{x}{1-r^2}, & \Gamma^y_{yy}=\frac{y}{1-r^2}
\end{pmatrix}
\ee
Recall that, in general,  the curvature (written with respect to a coordinate system $X_i=\frac{\partial}{\partial x_i}$) is defined as $\hat R(X_i,X_j)X_k:=\nabla_{X_j}\nabla_{X_i}X_k-\nabla_{X_i}\nabla_{X_j}X_k+\nabla_{[X_i,X_j]}X_k=\sum_l R^l_{ijk}X_l$, with $R^s_{ijk}=\sum_l\Gamma^l_{ik}\Gamma^s_{jl}-\sum_l\Gamma^l_{jk}\Gamma^s_{il}+\frac{\partial}{\partial x_j}\Gamma^s_{ik}-\frac{\partial}{\partial x_i}\Gamma^s_{jk}$. Moreover,  the sectional curvature of a two-dimensional subspace of the tangent space at a point which is spanned by $X,Y$ is given by $\hat K(X,Y)=\frac{g(R(X,Y)X,Y)}{|X|^2|Y|^2-g(X,Y)^2}$ and is independent of the choice of $X,Y$ (c.f. \cite{DoCarmo}, Chapter 4, Section 3, Proposition 3.1). In our case, letting $X=\frac{\partial}{\partial x}$ and $Y=\frac{\partial}{\partial y}$, $g_z(X,Y)=0$ (since the metric is diagonal), and $|X|^2=|Y|^2=\frac{1}{1-r^2}$. We compute 
{
\be{SU2sc}
\hat K(X,Y)=(1-r^2)^2g_z(\hat R(X,Y)X,Y)
\ee
$$=(1-r^2)^2(g_z(R^x_{xyx}\frac{\partial}{\partial x},\frac{\partial}{\partial y})+g_z(R^y_{xyx}\frac{\partial}{\partial y},\frac{\partial}{\partial y}))=
$$
$$(1-r^2)R^y_{xyx},$$
where
\be{Ryxyx}
R^y_{xyx}=
\ee
$$\Gamma^x_{xx}\Gamma^y_{yx}+\Gamma^y_{xx}\Gamma^y_{yy}-(\Gamma^x_{yx}\Gamma^y_{xx}+\Gamma^y_{yx}\Gamma^y_{xy})+\frac{\partial}{\partial y}\Gamma^y_{xx}-\frac{\partial}{\partial x}\Gamma^y_{yx}
$$
$$=\frac{-1}{(1-r^2)^2}(x^2-y^2-(-y^2+x^2)+1+y^2-x^2+1+x^2-y^2)=$$
$$\frac{-2}{(1-r^2)^2}$$
Therefore (\ref{SU2sc}) becomes 
\be{SU2sc2}
\hat K(X,Y)=(1-r^2)\frac{-2}{(1-r^2)^2}=\frac{-2}{1-r^2}
\ee
The open disc is two-dimensional, and so this shows that the the sectional curvature of $SU(2)_{\text{reg}}/K$ is nonpositive; hence Theorem \ref{negativo} shows that the cut locus of the identity must be contained in the inverse image (under the natural projection) of the boundary of the disc, i.e., the set of diagonal matrices in $SU(2)$. Furthermore from the calculation of the optimal trajectories for diagonal matrices in \cite{NoiAutomat}, it follows that the optimal trajectories cross the regular part. Therefore they lose optimality when they reach the set of diagonal matrices. In conclusion, in this case $K=SU(2)_{\text{sing}}$ coincides with the sub-Riemannian cut-locus. }


\noindent {\bf Acknowledgement} This research was partially  supported by NSF under Grant EECS-17890998

\end{document}